\documentclass[options]{article}
\usepackage{amsmath,amssymb}

\newcommand{\be}{\begin{equation}}
\newcommand{\ee}{\end{equation}}
\newcommand{\ba}{\begin{eqnarray}}
\newcommand{\ea}{\end{eqnarray}}
\newcommand{\baa}{\begin{eqnarray*}}
\newcommand{\eaa}{\end{eqnarray*}}
\newcommand{\bb}{\begin{thebibliography}}
\newcommand{\eb}{\end{thebibliography}}

\newcommand{\bi}[1]{\bibitem{#1}}
\newcommand{\lab}[1]{\label{#1}}
\newcommand{\re}[1]{(\ref{#1})}

\newcounter{my}
\newcommand{\he}%
   {\stepcounter{equation}\setcounter{my}%
   {\value{equation}}\setcounter{equation}0%
   }%
\newcommand{\she}%
   {\setcounter{equation}{\value{my}}%
    }%

\renewcommand\t{\tilde}

\newtheorem{pr}{Proposition}

\begin{document}

\vspace*{10mm}

\begin{center}

{\Large \bf Abstract "hypergeometric" orthogonal polynomials}

\vspace{5mm}

{\large \bf Alexei Zhedanov}

\medskip

{\em Donetsk Institute for Physics and Technology, Donetsk 83114,  Ukraine}

\end{center}


\begin{abstract}
We find all polynomials solutions $P_n(x)$ of the abstract "hypergeometric" equation
$L P_n(x) = \lambda_n P_n(x)$, where $L$ is a linear operator sending any polynomial of degree $n$ to a polynomial of the same degree with the property that $L$ is two-diagonal in the monomial basis, i.e. $L x^n = \lambda_n x^n + \mu_n x^{n-1}$ with arbitrary nonzero coefficients $\lambda_n, \mu_n$ . Under obvious nondegenerate conditions, the polynomial eigensolutions $L P_n(x) = \lambda_n P_n(x)$ are unique. The main result of the paper is a classification of all {\it orthogonal} polynomials  $P_n(x)$ of such type, i.e. $P_n(x)$ are assumed to be orthogonal with respect to a nondegenerate linear functional $\sigma$. We show that the only solutions are: Jacobi, Laguerre (correspondingly little $q$-Jacobi and little $q$-Laguerre and other special and degenerate cases), Bessel  and little -1 Jacobi polynomials.

\vspace{2cm}


\vspace{2cm}


\end{abstract}

\bigskip\bigskip


\newpage
\section{Introduction}
\setcounter{equation}{0} 
Let $L$ be an operator which transforms any polynomial of exact degree $n=0,1,2,\dots$ to a polynomial of the same exact degree $n$. Clealrly, this operator preserves any $n+1$-dimensional linear space of polynomials with degrees $\le n$ and hence there exists a system of eigenvalue monic polynomials $P_n(x)= x^n + O(x^{n-1})$ such that  
\be
L P_n(x) = \lambda_n P_n(x), \quad n=0,1,2,\dots \lab{eig_LP} \ee 
Under additional assumption that $\lambda_n \ne \lambda_m$ if $n \ne m$, the polynomals $P_n(x)$ are determined uniquely by equation \re{eig_LP}.

In oder to find explicit expression of the polynomials $P_n(x)$ let us assume that there exists a polynomial basis $\phi_n(x)=x^n + O(x^{n-1}), \: n=0,1,2,\dots$ such that the operator $L$ becomes low triangular and two-diagonal:
\be
L \phi_n(x) = \lambda_n \phi_n(x) + \mu_n \phi_{n-1}(x), \lab{hyp-2-diag-phi} \ee   
with some coefficients $\mu_n$.

If such basis is found then the eigenvalue problem \re{eig_LP} becomes simple. Indeed, we can write down the expansion
\be
P_n(x) = \sum_{s=0}^n A_{ns} \phi_s(x) \lab{P-hyp-phi} \ee
with unknown expansion coefficients $A_{ns}$. From \re{eig_LP} and \re{hyp-2-diag} we obtain 
\be
\frac{A_{n,s+1}}{A_{ns}} = \frac{\lambda_n-\lambda_s}{\mu_{s+1}}, \lab{rec_A} \ee
whence \be A_{ns}= A_{n,0} \:
\frac{(\lambda_n-\lambda_0)(\lambda_n-\lambda_1)\dots (\lambda_n
-\lambda_{s-1}) }{\mu_1 \mu_2 \dots \mu_s}, \quad s=1,2,\dots, n
\lab{expr_A} \ee The coefficient $A_{n0}$ can be chosen
arbitrarily. One possible choice is $A_{n0}=1$ for all $n$. This
corresponds to "hypergeometric-like" form of the polynomial
$P_n(x)$. Another choice is the monic form of the polynomial
$P_n(x)=x^n + O(x^{n-1})$. In this case \be A_{n0} = \frac{\mu_1
\mu_2 \dots \mu_n}{(\lambda_n-\lambda_0)(\lambda_n-\lambda_1)\dots
(\lambda_n -\lambda_{n-1}) } \lab{monic_A0} \ee  It is then 
convenient to present expansion coefficients in the form \be
A_{n,n-k} = \frac{\mu_n \mu_{n-1} \dots \mu_{n-k+1}}{(\lambda_n-
\lambda_{n-1})(\lambda_n- \lambda_{n-2}) \dots (\lambda_n-
\lambda_{n-k})}, \quad k=1,2,\dots, n \lab{A_monic} \ee
and we thus have the expression for the polynomial $P_n(x)$:
\be
P_n(x) = \phi_n(x) + \sum_{k=1}^n {\frac{\mu_n \mu_{n-1} \dots \mu_{n-k+1}}{(\lambda_n-
\lambda_{n-1})(\lambda_n- \lambda_{n-2}) \dots (\lambda_n-
\lambda_{n-k})} \phi_{n-k}(x)} \lab{P_n_mon} \ee
In the simplest case of the monomial basis $\phi_n(x) = x^n$ we have 
\be
Lx^n = \lambda_n x^n + \mu_n x^{n-1}, \lab{hyp-2-diag} \ee  
It is natural to call the operator $L$ the abstract "hypergeometric" operator.

Indeed, let us consider the classical hypergeomtric operator \cite{NSU}
\be 
L = x(1-x) \partial_x^2 + \left(\alpha+1 - (\alpha+\beta+2)x  \right)  \partial_x \lab{hyp_L_0} \ee
with arbitrary positive parameters $\alpha, \beta$. 

The property \re{hyp-2-diag} is almost obvious for this perator with 
\be
\lambda_n = -n( n+\alpha+\beta+1), \quad \mu_n = n(n+\alpha) . \lab{lm_hyp} \ee
From \re{expr_A} we obtain the explicit expression of (non-monic) $P_n(x)$ in terms of the Gauss hyergeometric function \cite{NSU}
\be
P_n(x) = \sum_{s=0}^n \frac{(-n)_s (n+\alpha+\beta+1)_s}{s! (\alpha+1)_s} x^s = {_2}F_1 \left( {-n, n+\alpha+\beta+1 \atop \alpha+1}; x \right) \lab{P-hyp} \ee
It appears that the eigenvalue polynomials $P_n(x)$ are orthogonal polynomials, i.e. they satisfy the orthogonality relation  
\be
\int_0^1 P_n(x) P_m(x) x^{\alpha}(1-x)^{\beta} dx =0, \quad n \ne m \lab{ort_Jac} \ee
In fact, $P_n(x)$ coincide with the classical Jacobi polynomials \cite{KLS}.

The q-hypergeometric differential operator \cite{KLS}, \cite{NSU} is another well known example of the operator $L$ with property \re{hyp-2-diag}. Define the operators
$$
T f(x) = f(qx), \: T^- f(x) = f(x/q), \; D_q f(x) =
\frac{f(xq)-f(x)}{x(q-1)}
$$
Then the second-order q-differential operator is \be L=
x(\gamma-x) T^- D_q^2 +(\alpha x + \beta) T^- D_q \lab{q-L_def}
\ee with 3 arbitrary parameters $\alpha,\beta, \gamma$. It is
easily verified that the operator $L$ has the property
\re{hyp-2-diag} with \be \lambda_n= q^{2-n}[n](q \alpha -[n-1]),
\; \mu_n= q^{2-n}[n](q \beta+ \gamma [n-1]), \lab{a_q_hyp} \ee
where $[n]=(q^n-1)/(q-1)$ is so-called "q-number". Similarly to the case of the ordinary hypergeometric operator, the operator \re{q-L_def} has orthogonal polynomials as eigenfunctions. These orthogonal polynomials are little q-Jacobi polynomials and their special and degenerate cases \cite{KLS}.

The main subject of the present paper is studying of abstract "hypergeometric" operators $L$ which act on monomials via formula \re{hyp-2-diag} with {\it a priori} arbitrary coefficients $\lambda_n, \mu_n$. The only restriction we are assuming is nondegeneracy of $\lambda_n$, i.e. we demand that 
\be
\lambda_n \ne \lambda_m, \quad \mbox{if} \quad n \ne m \lab{llne0} \ee
Under restriction \re{llne0} it is possible to find a unique set of monic polynomials $P_n(x) = x^n + O(x^{n-1})$ which are solutions of the problem \re{eig_LP}:
We have the explicit expression for the polynomials $P_n(x)$
\be
P_n(x) = x^n + \sum_{k=1}^n {\frac{\mu_n \mu_{n-1} \dots \mu_{n-k+1}}{(\lambda_n-
\lambda_{n-1})(\lambda_n- \lambda_{n-2}) \dots (\lambda_n-
\lambda_{n-k})} x^{n-k}} \lab{P_n_series} \ee

In what follows we will assume that $\lambda_0=\mu_0=0$. Indeed, condition $\mu_0=0$ is the truncation condition to prevent appearing of negative degrees of $x$. If $\lambda_0 \ne 0$ we can take the operator $L + const$ instead of the operator $L$ in order to achieve the condition $\lambda_0=0$.

The main purpose of the present paper is to find all possible operators $L$ such that corresponding eigenpolynomials $P_n(x)$ are orthogonal. The orthogonality condition (together with the natural condition of nondegeneracy) leads to very strong restrictions on the coefficients $\lambda_n, \mu_n$. We say that the sequences $\lambda_n$ and $\mu_n$ correspond to the admissible operator $L$ if the eigenpolynomials $P_n(x)$ are orthogonal and nondegenerate.

 The main result of our paper is that the only admissible operators $L$ 
are:

(i) $L$ is the q-hypergeometric operator of the second order. Corresponding polynomial solutions coincide with little q-Jacobi, q-Krawtchouk, little q-Laguerre, alternative q-Charlier and Stieltjes-Wiegert polynomials; 

(ii) $L$ is the ordinary hypergeometric operator (including the confluent case). This corresponds to the Jacobi, Laguerre and Bessel polynomials;

(iii) $L$ is of Dunkl type, i.e. $L$ is the first-order differential operator with the reflection operator $R$. This corresponds to the little -1 Jacobi polynomials.

The cases (i) and (ii) are well known \cite{KLS}. The orthogonal polynomials corresponding to the case (iii) were introduced and studied in \cite{VZ_little}.

The paper is organized as follows. In Section 2, we formulate and solve necessary and sufficient conditions for the coefficients $\lambda_n, \: \mu_n$ in order for polynomials $P_n(x)$  to be orthogonal. These conditions lead to a second-order linear recurrence relations for the coefficients $\lambda_n, \mu_n$. Solutions of these  relations strongly depend on the parameter $\Omega$. We distinguish cases when $\Omega \ne \pm 2$ (this is generic case) and $\Omega = \pm 2$. Solutions for these cases and corresponding families of orthogonal polynomials are considered in Sections 3--6. In Section 7, we analyze so-called "umbral classicality" of obtained families of orthogonal polynomials. In Section 8, we derive an equivalent necessary and sufficient condition for the operator to be admissible. This condition can be presented as an operator identity being a special case of the Askey-Wilson algebra AW(3).

\section{Necessary and sufficient conditions for orthogonality}
\setcounter{equation}{0}
In this section we derive and solve necessary and sufficient conditions for the coefficients $\lambda_n, \mu_n$ in order for polynomials $P_n(x)$ to be orthogonal and nondegenerate.

Recall that the polynomials $P_n(x)$ are orthogonal if there exists a nondegenerate linear functional $\sigma$ acting on the space of polynomials such that
\be
\langle \sigma, P_n(x) P_m(x) \rangle = 0, \quad n \ne m \lab{ort_P} \ee
where $\langle \sigma, \pi(x) \rangle$ stands for action of the functional $\sigma$ on a polynomial $\pi(x)$.

The functional $\sigma$ is completely determined by the moments
\be
\langle \sigma, x^n \rangle = c_n, \quad n=0,1,2,\dots \lab{mom_def} \ee
The linear functional
$\sigma$ is called nondegenerate if $\Delta_n \ne 0, \;
n=0,1,2,\dots$, where $\Delta_n = |c_{i+k}|_{i,k=0}^n$ are the Hankel determinants constructed from the moments. 

Conditions \re{mom_def} are equivalent to the conditions 
\be
\langle \sigma, P_n(x) \pi(x) \rangle = 0, \lab{ort_P2} \ee
where $\pi(x)$ is any polynomial of degree lesser than $n$.

Equivalently, the polynomials $P_n(x)$ are orthogonal iff they satisfy the three-term recurrence relation \cite{Chi}
\be
P_{n+1}(x) + b_n P_n(x) + u_{n} P_{n-1}(x) = xP_n(x), \quad n=1,2,\dots \lab{rec_P} \ee 
with some (complex, in general) coefficients $b_n, u_n$. The linear functional $\sigma$ is nondegenerate iff $u_n \ne 0, \; n=1,2,\dots$. In the special case when the coefficients $b_n$ are real and $u_n>0, \; n=1,2,\dots$ the polynomials are orthogonal with respect to a positive measure $d \mu(x)$ on the real axis \cite{Ismail}:
\be
\int_{a}^b P_n(x) P_m(x) d \mu(x) = h_n \delta_{nm}, \lab{ort_mu} \ee
where $h_0=1$ and $h_n=u_1 u_2 \dots u_n, \: n=1,2,\dots$. The integration limits $a$, $b$ in \re{ort_mu} may be either finite or infinite.

Assuming that the polynomials $P_n(x)$ are orthogonal (i.e. that they satisfy recurrence relation \re{rec_P}) we can find explicit expressions for the recurrence coefficients $b_n, u_n$ by comparing terms in front of monomials $x^n$ and $x^{n-1}$ in \re{rec_P} and by using \re{P_n_series}: 
\begin{subequations}\label{ub_expl}
\begin{align}
b_n & =  \frac{\mu_n}{\lambda_n - \lambda_{n-1}}  - \frac{\mu_{n+1}}{\lambda_{n+1} - \lambda_{n}} \label{b_n}\\
u_n & = -\frac{\mu_n b_n }{\lambda_n-\lambda_{n-1}} + \frac{\mu_n \mu_{n-1}}{(\lambda_n-\lambda_{n-1})(\lambda_n-\lambda_{n-2})} - \frac{\mu_n \mu_{n+1}}{(\lambda_{n+1}-\lambda_{n})(\lambda_{n+1}-\lambda_{n-1})}, \label{u_n}
\end{align}
\end{subequations}

Let $P_n(x), \: n=0,1,2,\dots$ be a set of monic orthogonal polynomials. Consider the polynomials  
\be
\tilde P_n(x) = \kappa^n \: P_n(x/\kappa) \lab{k_P} \ee 
It is clear that the polynomials $\tilde P_n(x)$ are monic $\tilde P_n(x) = x^n + O(x^{n-1})$. Moreover, these polynomials are orthogonal with respect to the functional $\tilde \sigma$ having the moments
\be
\tilde c_n = \kappa^n c_n \lab{t_cc} \ee
The polynomials $\t P_n(x)$ satisfy the three-term recurrence relation 
\be
\t P_{n+1}(x) + \t b_n \t P_n(x) + \t u_{n} \t P_{n-1}(x) = x \t P_n(x) \lab{t_rec} \ee
with
\be
\t b_n = \kappa b_n, \quad \t u_n = \kappa^2 u_n \lab{t_bu} \ee
Although formulas \re{t_cc}, \re{t_bu} are almost obvious, they can be useful in reducing given orthogonal polynomials to a more simple form.

Assume that $L$ is any operator which preserves degree of polynomials. Then it is possible to construct a unique system of monic eigenpolynomials $P_n(x)$ satisfying \re{eig_LP} (of course, it is assumed that the spectrum $\lambda_n$ is nondegenerate). In \cite{VYZ} it was shown that the polynomials $P_n(x)$ will be orthogonal with respect to a linear functional $\sigma$ if and only if the condition 
\be
\langle \sigma, g(x) L f(x) \rangle =  \langle \sigma, f(x) L g(x) \rangle \lab{gLf} \ee
holds for any pair of polynomials $f(x), g(x)$. 
Condition \re{gLf} has a simple meaning that the operator $L$ is symmetric on the space of polynomials with respect to the functional $\sigma$. See also \cite{Duran} where the same condition is derived in case if $L$ is a higher-order difference operator (in fact, this condition is valid for {\it any} linear operator $L$ preserving polynomiality such that $\deg(L p(x)) \le \deg(p(x))$ for any polynomial $p(x)$).

Equivalently, condition \re{gLf} can be presented in the form of an infinite set of conditions
\be
\langle \sigma, x^m L x^n \rangle =  \langle \sigma, x^n L x^m \rangle \lab{mLn} \ee
which should be valid for all possible nonnegative integers $m,n=0,1,2,\dots$

Taking into account relation \re{hyp-2-diag} and remembering definition \re{mom_def} of moments $c_n$ we can present conditions \re{mLn} in the form
\be
(\lambda_n-\lambda_m)c_{n+m} + (\mu_n-\mu_m)c_{n+m-1}=0, \quad m,n=0,1,2,\dots
\lab{ln_cond} \ee

Relations \re{ln_cond} are necessary and sufficient conditions for orthogonality of eigenpolynomials $P_n(x)$. These conditions contain 3 unknown sequences: $\lambda_n, \: \mu_n, \: c_n, \; n=0,1,2,\dots$ with initial conditions
\be
\lambda_0=\mu_0=0, \; c_0=1 \lab{ini_lmc} \ee
First of all we can conclude from \re{ln_cond} that $\mu_n \ne \mu_m$ if $n \ne m$.
Indeed, assume that $\mu_j=\mu_k$ for a pair of nonnegative
integers $j,k$ such that $j>k \ge 0$. Then from \re{ln_cond} it
follows that $c_{j+k}=0$ (because $\lambda_j \ne \lambda_k$ by our
assumption). By induction from the same condition we then find
that $c_{n} =0, n=j+k, j+k+1, j+k+2, \dots$. But this contradicts
to nondegenerate condition $\Delta_n \ne 0$ for all $n>0$.

Hence we have that both coefficients $\lambda_n$ and $\mu_n$ are
nondegenerate: \be \lambda_n \ne \lambda_m, \quad \mu_n \ne \mu_m
\quad \mbox{if} \quad n \ne m \lab{lanu_ndeg} \ee

From this condition it follows immediately that that all moments
are nonzero: $c_n \ne 0, \; n=0,1,2,\dots$. Indeed, if one assumes
that $c_j=0$ for some $j>0$ then from \re{ln_cond} we find that
$c_n=0$ for all $n \ge j$ which contradicts to nondegenerate
condition $\Delta_n \ne 0$.

Hence we have the condition \be \frac{\lambda_n -
\lambda_m}{\mu_n-\mu_m} = g_{n+m}, \lab{I_I_3} \ee where
$$
g_n = -c_{n-1}/c_n.
$$

Note also that $\mu_n$ cannot be a linear function of $\lambda_n$, i.e.
condition $\mu_n= \alpha \lambda_n + \beta, \; n=1,2,3,\dots$ is forbidden.
Otherwise from \re{I_I_3} it follows $g_n \equiv const$ which
leads to a degeneration $\Delta_n =0$.

It is easy to see from conditions \re{I_I_3} that 5 parameters $\lambda_1, \mu_1, \lambda_2, \mu_2, \lambda_3$ can be chosen arbitrarily. For the parameter $\mu_3$ we have from \re{I_I_3}
\be
\mu_3=\lambda_3 (\mu_2-\mu_1)/(\lambda_2-\lambda_1) \lab{nu3} \ee 
Then all further parameters $\lambda_4, \mu_4, \lambda_5, \mu_5, \dots$ are determined uniquely from conditions \re{I_I_3}. We shall thus assume that these 5 parameters are fixed as arbitrary parameters (with obvious restrictions of nondegeneracy).

Eliminating $g_{m+n}$ from \re{I_I_3}  we arrive at the relations \be \frac{\lambda_{n+1} -
\lambda_k}{\mu_{n+1} - \mu_k} = \frac{\lambda_{n} -
\lambda_{k+1}}{\mu_{n} - \mu_{k+1}} \lab{ln_nk} \ee which are valid
for $n,k=0,1,2,\dots$. It is easily seen that relations \re{ln_nk} are equivalent to \re{I_I_3}, so we can restrict ourselves only with relations \re{ln_nk} which are more convenient because they involve only two unknown sequences $\lambda_n, \mu_n$.

Putting $k=0,1,2$ and taking into account initial conditions
\re{ini_lmc} we have
\be
\frac{\lambda_{n+1}}{\mu_{n+1}} = \frac{\lambda_{n} - \lambda_1}{\mu_{n} -
\mu_1}, \quad \frac{\lambda_{n+1} - \lambda_1}{\mu_{n+1} - \mu_1} =
\frac{\lambda_{n} - \lambda_{2}}{\mu_{n} - \mu_{2}}, \quad  \frac{\lambda_{n+1} - \lambda_2}{\mu_{n+1} - \mu_2} =
\frac{\lambda_{n} - \lambda_{3}}{\mu_{n} - \mu_{3}}, \; 
 \lab{eq_ln} \ee
Eliminating variables $\mu_n$ and $\mu_{n+1}$ from 3 equations \re{eq_ln}, we obtain the equation
\be
A_1 \lambda_n^2 + A_2 \lambda_{n+1}^2 + A_3 \lambda_n \lambda_{n+1} + A_4 \lambda_n + A_5 \lambda_{n+1} + A_6=0, \quad n=4,5,6,\dots  \lab{qua_l} \ee
where
$$
A_1= \mu_2 \lambda_1-\mu_1 \lambda_2, \quad A_2 = A_1  + \lambda_3(\mu_1-\mu_2) + \mu_3 (\lambda_2-\lambda_1)
$$ 
$$
A_3= \mu_1 (\lambda_3-\lambda_2) + \lambda_1 (\mu_2 -\mu_3) 
$$
$$
A_4=(\lambda_1+\lambda_2)(\lambda_2 \mu_1 - \lambda_1 \mu_2) + \lambda_1 (\lambda_2 \mu_3 - \lambda_3 \mu_2), \; A_5 = A_4 -(\lambda_1 + \lambda_2)(\lambda_3 (\mu_1-\mu_2) + \mu_3 (\lambda_2-\lambda_1))
$$
$$
A_6= \lambda_1 \left( \lambda_1 (\lambda_3 \mu_2 - \lambda_2 \mu_3) + \lambda_2 (\lambda_1 \mu_2 - \lambda_2 \mu_1) \right)
$$
We notice that $A_1 \ne 0$. Indeed, if one assumes that $A_1=0$ then by induction it is possible to show that $\lambda_n /\mu_n = \lambda_1/\mu_1, \: n=2,3,4,\dots$ which leads to a degeneration.

Moreover, from \re{nu3} we have that $A_2=A_1$ and $A_5=A_4$.

Hence one can rewrite equation \re{qua_l} in the form
\be
\lambda_n^2 + \lambda_{n+1}^2 + B_1 \lambda_n \lambda_{n+1}  + B_2 (\lambda_n + \lambda_{n+1}) + B_3 =0 \lab{qua_sym_l} \ee
where
\be
B_1= 1+ \frac{\mu_1 \lambda_3 - \mu_3 \lambda_1}{\mu_2 \lambda_1 - \mu_1 \lambda_2} = 1+ \frac{\lambda_3}{\lambda_1-\lambda_2} \lab{B_1} \ee
 (the last equality in \re{B_1} follows from  \re{nu3}).

 Using the same relation \re{nu3} we can present the coefficients $B_2, B_3$ in the form
\be
B_2=\frac{\lambda_1^2 -\lambda_2^2 + \lambda_1 \lambda_3}{\lambda_2-\lambda_1}  \lab{B_2} \ee
and 
\be
B_3 = \lambda_1 \lambda_2 + \frac{\lambda_1^2 \lambda_3}{\lambda_1-\lambda_2} \lab{B_3} \ee 
Due to obvious symmetry between $\lambda_n$ and $\mu_n$ one can obtain similar equation for $\mu_n$:
\be
\mu_n^2 + \mu_{n+1}^2 + C_1 \mu_n \mu_{n+1}  + C_2 (\mu_n + \mu_{n+1}) + C_3 =0 \lab{qua_sym_nu} \ee 
where the coefficients $C_1,C_2,C_3$ have expressions \re{B_1}- \re{B_3} with $\lambda_i$ replaced with $\mu_i$. It is important to note that coefficients $B_1$ and $C_1$ coincide
\be
C_1=B_1 =  -\Omega = 1+ \frac{\lambda_3}{\lambda_1-\lambda_2} \lab{CB_1} \ee
Quadratic relations \re{qua_sym_l} and \re{qua_sym_l} contain 5 free parameters $\Omega, B_2, B_3, C_2, C_3$. Putting $n=3$ we find from quadratic equation \re{qua_sym_l} that there are 2 possible solutions for $\lambda_4$. But due to obvious symmetry between $\lambda_n$ and $\lambda_{n+1}$ in \re{qua_sym_l}, one of these solutions corresponds to $\lambda_2$. This solution should be excluded because $\lambda_4=\lambda_2$ means the degeneration. We thus have only one solution for $\lambda_4$. By induction, it is easy to show that all values $\lambda_n, \: \mu_n, \: n=4,5,6,\dots$ are determined uniquely.

From \re{qua_sym_l} and \re{qua_sym_nu} (under conditions $\lambda_{n+1} \ne\lambda_{n-1}$ and $\mu_{n+1} \ne\mu_{n-1}$) we derive the linear recurrence relations
\be
\lambda_{n+1}+\lambda_{n-1} -\Omega  \lambda_n +B_2 =0, \quad n=1,2,\dots \lab{lin_lam} \ee  
and
\be
\mu_{n+1}+\mu_{n-1} -\Omega  \mu_n +C_2 =0, \quad n=1,2,\dots \lab{lin_nu} \ee  
Recurrence equations \re{lin_lam} and \re{lin_nu} are well known. They arise, e.g. in
describing so-called "Askey-Wilson grids" \cite{AW}, \cite{NSU},
\cite{VZ_Bochner}. 

On the other hand, conditions \re{lin_lam} and \re{lin_nu} are not only necessary but also {\it sufficient} 
for our problem. This will be seen from the results of the next four  sections where we classify all possible solutions of the recurrence relations  \re{lin_lam} and \re{lin_nu}.

We thus can formulate the main result of this section:
\begin{pr}
The abstract "hypergeomtric" operator $L$ defined by \re{hyp-2-diag} has orthogonal polynomials $P_n(x)$ as eigensolutions if and only if the coefficients $\lambda_n$ and $\mu_n$ are solutions of equations \re{lin_lam} and \re{lin_nu} with arbitrary parameters $B_2,C_2$. Additional restriction for the nondegeneracy is: the linear relation  $\alpha \lambda_n + \beta \mu_n+ \gamma =0$ is forbidden (in case if at least on of the constants $\alpha$ and $\beta$ is nonzero). 
\end{pr}

\section{Classification of admissible solutions. $\Omega>2$}
\setcounter{equation}{0}
General solution of recurrence relations  \re{lin_lam} and \re{lin_nu} depends on the parameter $\Omega$.

Assume first that $\Omega>2$. In this case we can put 
\be
\Omega = q+ q^{-1}, \lab{Om_q} \ee
where $0<q<1$ is a real positive parameter. 

We can take $\lambda_1, \lambda_2, \mu_1, \mu_2$ as independent initial parameters. Then 
\be
\lambda_3= \frac{(\lambda_2-\lambda_1)(1-q^3)}{q(1-q)}, \quad  \mu_3= \frac{(\mu_2-\mu_1)(1-q^3)}{q(1-q)} \lab{ln_3} \ee

Generic solution with initial conditions $\lambda_0=\mu_0=0$ in this case is 
\be
\lambda_n = L_1 (q^n -1) + L_2 (q^{-n}-1) \lab{lam_n} \ee
and
\be
\mu_n = M_1 (q^n -1) + M_2 (q^{-n}-1), \lab{nu_n} \ee
where $L_1, L_2$ depend on initial conditions $\lambda_1, \lambda_2$ (correspondingly, $M_1, M_2$ depend on initial conditions $\mu_1, \mu_2$):
\be
L_1= \frac{q \lambda_2 -(q+1) \lambda_1}{(q+1)(q-1)^2}, \quad  L_2= \frac{q^2(\lambda_2 -(q+1) \lambda_1)}{(q+1)(q-1)^2} \lab{L_12} \ee
and
\be
M_1= \frac{q \mu_2 -(q+1) \mu_1}{(q+1)(q-1)^2}, \quad  M_2= \frac{q^2(\mu_2 -(q+1) \mu_1)}{(q+1)(q-1)^2} \lab{M_12} \ee
The parameters $B_2,C_2$ in equations \re{lin_lam} and \re{lin_nu} are 
\be
B_2=-q^{-1} (q-1)^2 (L_1+L_2), \; C_2=-q^{-1} (q-1)^2 (M_1+M_2) \lab{BC_q} \ee 
It is easily seen that condition \re{I_I_3} is fulfilled with 
\be
g_n = -\frac{c_{n-1}}{c_n} = \frac{L_2-L_1q^n}{M_2-M_1q^n} \lab{g_n_q} \ee
Generic solution depends on 4 arbitrary constants $L_1,L_2,M_1,M_2$. In fact, only two constants can be considered as independent parameters because rescaling transformations $\lambda_n \to \kappa_1 \lambda_n, \: \mu_n \to \kappa_2 \mu_n$ lead to a rescaling of the argument of the orthogonal polynomials $P_n(x)$.

\vspace{5mm}

(i)   Assume that $L_1L_2M_1M_2 \ne 0$. We can fix two parameters as, e.g. $L_2=1, \: M_2=-1$. It is convenient to introduce two independent parameters $a,b$ such that $L_1=abq, \: M_1=-a$. 

Corresponding operator $L$ can be presented as the difference q-hypergeometric operator
\be
L f(x) = a(bq-x^{-1})(f(qx)-f(x)) + (1-x^{-1})(f(x/q)-f(x)) \lab{little_L} \ee
Using formulas \re{ub_expl} we find explicit expressions for the recurrence coefficients
\be
u_n = A_{n-1} C_n, \quad b_n = A_n + C_n, \lab{ub_little} \ee 
where
\be
A_n = q^n \: \frac{(1-a q^{n+1})(1-ab q^{n+1})}{(1-abq^{2n+1})(1-abq^{2n+2})}, \quad C_n = aq^n \: \frac{(1-q^{n})(1-b q^{n})}{(1-abq^{2n+1})(1-abq^{2n})} \lab{AC_little} \ee
Formulas \re{ub_little} coincide with the recurrence coefficients for the little q-Jacobi polynomials \cite{KLS}. 

For $0<a<q^{-1}, \; b<q^{-1}$ these polynomials are orthogonal on the infinite set of the points
\be
x_s=q^s, \quad s=0,1,2,\dots \lab{x_s_LJ} \ee
with the weights
\be
w_s = (aq)^s \: \frac{(bq;q)_s}{(q;q)_s}, \lab{w_LJ} \ee
where $(a;q)_s =(1-a)(1-aq) \dots (1-aq^{s-1})$ is the shifted q-factorial (q-Pochhammer symbol).

If $b= q^{-N-1}$ then a finite set of orthogonal polynomials appears. They are orthogonal on the set of $N+1$ points of the real axis
\be
\sum_{s=0}^N  {(aq)^s \: \frac{(q^{-N};q)_s}{(q;q)_s} P_n(q^s)P_m(q^s) } =0, \quad n \ne m \lab{finite_lJ} \ee
These polynomials can be identified with q-Krawtchouk polynomials \cite{KLS}. For positivity of the weights in \re{finite_lJ} it is necessary that $a<0$. Note that these q-Krawtchouk polynomials differ from the "standard" q-Krawtchouk polynomials (see, e.g. \cite{KLS}) by the change $q \to q^{-1}$.

We thus see that the generic case of the solutions for $\lambda_n, \mu_n$ corresponds to the little q-Jacobi polynomials.

(ii) Assume that $L_1=0$ and $M_1 M_2 \ne 0$. Then the only one parameter is essential. Without loss of generality we can put $L_2=-1, M_1=a, M_2=1$, where $a$ is the only essential parameter. From formulas \re{ub_expl} we obtain
\be
u_n= aq^{2n-1} (1-q^n)(1-aq^n), \quad b_n = (1+a)q^n - a(1+q) q^{2n} \lab{ub_l_Lag} \ee
These formulas coincide with expressions for the recurrence coefficients of the little q-Laguerre polynomials \cite{KLS}.

(iii) Assume that $L_2=0$ and $M_1 M_2 \ne 0$. Again there is only one essential parameter and we can put $M_2=1, \: L_1 =M_1=a$. For the recurrence coefficients we obtain the expressions
\be
u_n = \frac{q^{-4n+1}}{a^2} (1-q^n ) (1-a q^n), \quad b_n = \frac{q+1}{a}q^{-2n+1} -\frac{a+1}{a} q^{-n} \lab{rec_q-L} \ee
These correspond to the recurrence coefficents of the q-Laguerre polynomials \cite{KLS}.

(iv) Assume that $M_1=0$ and $L_1 L_2 \ne 0$. We can put $L_1=a, \: L_2=-1, \: M_2=1$ with the only essential parameter $a$.
The recurrence coefficients are 
\begin{subequations}\label{ub_alter}
\begin{align}
u_n = \frac{aq^{3n-2} (1-q^n)(1+a q^{n-1})}{(1+aq^{2n})(1+aq^{2n-2})(1+aq^{2n-1})^2}, \\ 
b_n =\frac{1-q^{n+1}}{(1-q)(1+aq^{2n+1})} -  \frac{1-q^{n}}{(1-q)(1+aq^{2n-1})}
\end{align}
\end{subequations}
These coefficients correspond to alternative q-Charlier polynomials.

(v) Assume that $L_1=M_2=0$ and $L_2M_1 \ne 0$. Without loss of generality we can put $L_2=M_1=1$ whence
\be
\lambda_n = q^{-n}-1, \quad \mu_n = q^n -1 \lab{lm_SW} \ee
The recurence coefficients 
\be
u_n= q^{1-4n} (1-q^n) , \quad b_n = (q+1) q^{-2n-1} - q^{-n} \lab{ub_SW} \ee
correspond to the Stieltjes-Wigert polynomials \cite{KLS}.

The case $L_2=M_1=0$ and $L_1 M_2 \ne 0$ is "dual" with respect to the previous case: it corresponds to the transformation $q \to q^{-1}$.

Remaining cases, say $L_1=M_1=0$ or $M_1=M_2=0$ correspond to degenerate polynomials because in these cases there exists the linear dependence $\mu_n=\alpha \lambda_n + \beta$ which leads to a degeneration. This can be confirmed by direct calculation of the recurrence coeficient $u_n$ using formulas \re{ub_expl}. We obtain $u_n =0$ for all $n=1,2,\dots$ which coresponds to a degenerate case.

\section{Admissible solutions. $\Omega<-2$ and $-2<\Omega<2$}
\setcounter{equation}{0}
The case $\Omega<-2$ is very close to the already considered case $\Omega>2$. All formulas of the previous section remain valid if one changes $ q \to -q$. This is equivalent to the change $\Omega \to -\Omega$. However, the spectral properties of the corresponding orthogonal polynomials will be slightly different. For example, for the case (i) we have the little q-Jacobi polynomials with $-1<q<0$. In this case the spectrum will formally be the same $x_s=q^s, \: s=0,1,2,\dots$ but now this set is a union of two geometric series in the interval $[-1,1]$ with the concentrated point $x_{\infty}=0$. For the positivity of the weight function it is necessary that $a<0$ and $|a|<|q|^{-1}$.

The case $-2<\Omega<2$ corresponds to trigonometric expressions of the recurrence coefficients $u_n, b_n$.

Indeed, in this case we can put 
\be
q=e^{2i\omega} \lab{q_trig} \ee
with some real parameter $\omega$.

If we demand that the coefficients $\lambda_n, \mu_n$ are real then from \re{lam_n}, \re{nu_n} it is clear that necessarily $L_2=L_1^*, \: M_2 = M_1^*$. This is equivalent to the representation 
\be
\lambda_n = \sin \omega n \sin (\omega(n+\alpha+\beta+1)), \quad  \mu_n = -\sin \omega n \sin (\omega(n+\alpha)) \lab{la_mu_trig} \ee
with some real parameters $\alpha,\beta$.

In turn, substituting these formulas into \re{ub_expl} we obtain
\be
b_n =  
\frac{\sin \omega (n+1)  \sin \omega(n+\alpha+1)} {\sin \omega \sin \omega (2n+\alpha+\beta+2)} -\frac{\sin \omega n  \sin \omega(n+\alpha)} {\sin \omega \sin \omega (2n+\alpha+\beta)}   \lab{b_trig} \ee 
and 
\be
u_n={\frac {\sin \omega\,n  \sin \omega\, \left( n+
\beta \right)   \sin \omega\, \left( n+\alpha \right) 
  \sin  \omega\, \left( n+\alpha+\beta \right) 
  }{\sin  \omega\, \left( 2\,n+\alpha+\beta+1 \right)   
\sin  \omega\, \left( 2\,n+\alpha+\beta-1 \right)    
\sin^2  \omega\, \left( 2\,n+\alpha+\beta   \right)  }}
\lab{u_trig} \ee
It is impossible to provide the positivity condition $u_n>0$ for all $n$. This means that there is no a  positive measure on the real line such that the polynomials $P_n(x)$ are orthogonal with respect to it.

There is a simple special case when $\beta=-N-1$. In this case $u_{N+1}=0$ and we have a finite set of polynomials $P_n(x)$orthogonal on the vertices of the regular  $N+1$-gon on the unit circle    
\be
\sum_{s=0}^N w_s P_n(x_s) P_m(x_s) =0, \quad n \ne m \lab{trig_ort} \ee
where 
\be
x_s = \exp(i \omega (-N + 2is)), \quad w_s = \frac{(q^{-N};q)_s q^{(\alpha+1)s}}{(q;q)_s} \lab{xw_trig} \ee
As expected, the weights $w_s$ are not positive (or even real) parameters. 

Another interesting case of a finite orthogonality appears when $q$ is a root of unity. For example, one can take $\omega = \frac{\pi}{N+1}$. In this case $u_{N+1}=0$ and we again have a finite set of the orthogonal polynomials which are orthogonal on the circle (in general, with non-unit radius). In this case again there is no a positivity property. We will not consider this case in details.   

General discussion concerning the Askey-Wilson polynomials for $q$ a root of unity can be found in \cite{SZ_root}.

\section{Admissible solutions. $\Omega=2$}
\setcounter{equation}{0}
Consider the case $\Omega=2$. It is easily seen that generic solutions of equations  \re{lin_lam} and \re{lin_nu} with initial conditions $\lambda_0=\mu_0=0$ are
\be
\lambda_n = L_2 n^2 + L_1 n, \quad \mu_n = M_2 n^2 + M_1 n \lab{lm_Om=2} \ee 
The parameters $B_2,C_2$ in equations \re{lin_lam} and \re{lin_nu} are 
\be
B_2=-2L_2, \; C_2=-2M_2 \lab{BC_1} \ee
It is clear that necessary and sufficient condition \re{I_I_3} is valid leading to the equation for the moments 
\be
\frac{c_n}{c_{n-1}} = - \frac{M_2n + M_1}{L_2 n + L_1} \lab{mom_Jac} \ee
Depending on the choice of the parameters $L_{1,2}, M_{1,2}$ we can distinguish 3 possibilities.

(i) Assume that $L_2 M_2 \ne 0$. Then it is possible to choose $M_2=2, \: L_2=1, \: M_1=2 \alpha, \: L_1=\alpha+\beta+1$. 
This correspods to the hypergeometric operator
\be
L= x(1-x)\partial_x^2 +\left( \alpha+1 - (\alpha+\beta+2)x  \right) \partial_x \lab{hyp_L} \ee
The recurrence coefficients can be computed through formulas \re{ub_expl}
\begin{subequations}\label{ub_Jac}
\begin{align}
u_n = \frac{n(n+\alpha)(n+\beta)(n+\alpha+\beta)}{(2n+\alpha+\beta+1)(2n+\alpha+\beta-1)(2n+\alpha+\beta)^2}, \\ 
b_n = \frac{\alpha^2-\beta^2}{2(2n+\alpha+\beta)(2n+\alpha+\beta+2)} +1/2
\end{align}
\end{subequations}
These recurrence coefficients correspond to the Jacobi polynomials $P_n^{(\alpha,\beta)}\left(\frac{x+1}{2} \right)$ \cite{KLS} . They are orthogonal on the interval $[0,1]$ with respect to the weight function
\be
w(x) = x^{\alpha} (1-x)^{\beta} .\lab{w_J} \ee

(ii) Assume that $L_2=0$ and $M_2 \ne 0$. There is only one essential parameter, say $\alpha$, and one can choose $L_1=1, \: M_2=-1, \: M_1=-\alpha$. 
Corresponding operator $L$ is second-order differential operator 
\be
L=-x \partial_x^2  +(x-\alpha-1) \partial_x \lab{L_Lag} \ee 
The recurrence coefficients
\be
u_n = n(n+\alpha), \quad b_n = 2n+\alpha+1 \lab{ub_Lag} \ee 
correspond to the Laguerre polynomials \cite{KLS}. They are orthogonal on the semi-axis $[0,\infty]$ with respect to the weight
\be
w(x) = x^{\alpha} e^{-x} \lab{w_Lag} \ee

(iii) Assume that $M_2=0$ and $L_2 \ne 0$. One can put 
\be
\lambda_n = n(n+a-1), \quad \mu_n =2 \lab{lm_Bes} \ee
with the only real parameter $a$.

Corresponding operator $L$ is second-order differential operator 
\be
L = x^2 \partial_x^2 +(ax+2) \partial_x \lab{L_Bessel} \ee

The recurrence coefficients
\begin{subequations}\label{ub_Bes}
\begin{align}
u_n &= -\frac{4n(n+a-2)}{(2n+a-1)(2n+a-3)(2n+a-2)^2}, \\ 
b_n &= \frac{a-2}{2n+a} - \frac{a-2}{2n+a-2}
\end{align}
\end{subequations}
correspond to the generalized Bessel polynomials \cite{Chi}. Note that $u_n$ cannot be positive for all $n$. This means that the Bessel polynomials cannot be orthogonal with respect to a positive measure on the real axis \cite{Chi}.

\section{Admissible solutions. $\Omega=-2$}
\setcounter{equation}{0}
If $\Omega=-2$ then generic solutions of equations  \re{lin_lam} and \re{lin_nu} with initial conditions $\lambda_0=\mu_0=0$ are
\be
\lambda_n = (-1)^n (L_1 n + L_0) - L_0, \quad  \mu_n = (-1)^n (M_1 n + M_0) - M_0 \lab{lm_-1} \ee
with four arbitrary parameters $L_1,L_0,M_1,M_0$.
The parameters $B_2,C_2$ in equations \re{lin_lam} and \re{lin_nu} are 
\be
B_2=4L_0, \; C_2=4M_0 \lab{BC--1} \ee
Necessary and sufficient condition \re{I_I_3} is valid leading to the equation for the moments 
\be
\frac{c_n}{c_{n-1}} = - \frac{(-1)^n(M_1 n + M_0) - M_0}{(-1)^n(L_1 n + L_0) - L_0} \lab{mom-1} \ee

Note that necessarily $L_1 M_1 \ne 0$ because otherwise either $\lambda_{2n+1}=0$ or $\mu_{2n+1}=0$ which is forbidden by nondegeneracy conditions.

Hence two parameters $L_1,M_1$ can be chosen as fixed nonzero constants, say $L_1=-2, \: M_1=2$. The two remaining parameters can be parametrized as $L_0 = -\alpha-\beta-1, \: M_0=\alpha$ with two arbitrary real parameters $\alpha,\beta$.

Then it is easy to see that corresponding operator $L$ has the expression
\be
L=2(1-x)\partial_x R +(\alpha+\beta+1-\alpha x^{-1})(1-R), \lab{L-1} \ee
where $R$ is the reflection operator, i.e. $R f(x) = f(-x)$.

As was shown in \cite{VZ_little} the operator \re{L-1} is the Dunkl type differntial operator such that its polynomial
eigenfunctions $P_n(x)$
\be
L P_n(x) = \lambda_n P_n(x) \lab{LP_-1} \ee
coincide with the little -1 Jacobi polynomials.

Using formulas \re{ub_expl} we can calculate the recurrence coefficients
\be u_{n} = \frac{(n+(1-\theta_n)\alpha)(n+\beta+
\theta_n \alpha)}{(2n+\alpha+\beta)^2}, \quad b_n = (-1)^n \:
\frac{(2n+1)\alpha + \alpha \beta +\alpha^2 +(-1)^n
\beta}{(2n+\alpha+\beta)(2n+2+\alpha+\beta)}, \lab{lqJ-1_ub} \ee
where
$$
\theta_n= \frac{1+(-1)^n}{2}
$$
is the characteristic function of even numbers. As expected, these coefficients coincide with recurrence coefficients of the little -1 Jacobi polynomials \cite{VZ_little}.

The little -1 Jacobi polynomials are orthogonal on the interval $[-1,1]$ with respect to the weight function \cite{VZ_little}
\be
w(x) = |x|^{\alpha}(1+x) (1-x^2)^{\frac{\beta-1}{2}} \lab{w_little} \ee
There is a special case corresponding to $M_0=0$, i.e. $\alpha=0$. This means that $\mu_n = 2(-1)^n \: n$. As was shown in \cite{VZ_little} this case corresponds to the ordinary Jacobi polynomials $P_n^{(a,a+1)}(x)$ with $a=(\beta-1)/2$. It is interesting to note that the classical Jacobi polynomials $P_n^{(a,a+1)}(x)$ satisfy "unusual" eigenvalue problem \re{LP_-1} with the Dunkl type operator \re{L-1}. Explanation of this phenomenon can be found in \cite{VZ_little}.

\section{The umbral "classical" polynomials}
\setcounter{equation}{0}
There is an interesting relation of the above approach with the so-called "umbral calculus" \cite{Roman}. One of the most important object in the umbral calculus is the formal derivative operator $\cal D$ which is defined on the space of polynomials by its action on monomials
\be
{\cal D} x^n = d_n x^{n-1}, \quad n=0,1,2,\dots\lab{def_D} \ee
where $d_n$ is an arbitrary sequence of complex numbers with the restrictions $d_0=0$ and $d_n \ne 0, \: n=1,2,\dots$  Clearly, the operator $\cal D$ decreases the degree of any polynomial by one. The obvious example of the formal derivative operator is the ordinary derivative operator ${\cal D} f(x)  = f'(x)$ (in this case $d_n =n$). Another simple example is the q-derivative operator
\be
D_{q} f(x) =\frac{f(xq)-f(x)}{x(q-1)}. \lab{D_q} \ee
In this case $d_n = (q^n-1)/(q-1)$.

We say that the set of monic orthogonal polynomials $P_n(x), \: n=0,1,2,\dots $ satisfies the umbral classical property if the new set of monic polynomials
\be
\tilde P_n(x) = \frac{{\cal D} P_{n+1}(x)}{d_{n+1}}, \quad n=0,1,2,\dots \lab{tPD} \ee
is another set of orthogonal polynomials (see e.g. \cite{Khol} for further details).

The ordinary classical polynomials  (i.e. the Jacobi, Laguerre, Bessel and Hermite polynomials) satisfy this definition with ${\cal D} = \partial_x$. This statement is known as the Hahn theorem \cite{Al-Salam}.

It is easy to see that the orthogonal polynomials considered in this paper do satisfy the umbral classical property.

Indeed, consider the explicit form of series expansion of orthogonal polynomials   \re{A_monic}. Applying the operator $\cal D$ to $P_{n+1}(x)$ we get
\be \tilde P_n(x) = x^n + \sum_{k=1}^n {\frac{\mu_{n+1} \mu_{n} \dots \mu_{n-k+2}}{(\lambda_{n+1}-
\lambda_{n})(\lambda_{n+1}- \lambda_{n-1}) \dots (\lambda_{n+1}-
\lambda_{n-k+1})} \frac{d_{n+1-k}}{d_{n+1}} x^{n-k}} \lab{D_P} \ee
In all cases corresponding to $\Omega \ne \pm 2, \: \Omega=2, \: \Omega=-2$ it is possible to rewrite \re{D_P} in the equivalent form
\be 
\tilde P_n(x) = x^n + \sum_{k=1}^n {\frac{\tilde  \mu_n \tilde \mu_{n-1} \dots \tilde \mu_{n-k+1}}{(\tilde \lambda_n-
\tilde \lambda_{n-1})(\tilde \lambda_n- \tilde \lambda_{n-2}) \dots (\tilde \lambda_n-
\tilde \lambda_{n-k})} x^{n-k}} \lab{equiv_Pt} \ee
with new coefficients $\tilde \lambda_n, \: \tilde \mu_n $ belonging to the same admissible class of solutions of equation \re{I_I_3}. This means that the new polynomials $\tilde P_n(x)$ are orthogonal and satisfy similar eigenvalue equation 
\be
\tilde L \tilde P_n(x) = \tilde \lambda_n \tilde P_n(x), \lab{tLP} \ee
where the operator $\tilde L$ is defined as
\be
\tilde L x^n = \tilde \lambda_n x^n + \tilde \mu_n x^{n-1} \lab{tL_x} \ee
The concrete form of the operator $\cal D$ depends on the value of $\Omega$.

In more details, when  $\Omega=2$, we have $d_n=n$, i.e. $\cal D$ coincides with the ordinary derivative operator. When $\Omega \ne \pm 2$, we have $d_n = (q^n-1)/(q-1)$, i.e. in this case $\cal D$ is the q-derivative operator (see \cite{Al-Salam} for details).

Finally when $\Omega=-2$ we have 
\be
d_n = n + \nu \: (1-(-1)^n) \lab{d_nu} \ee
with an appropriate parameter $\nu$. In this case $\cal D$ is the classical Dunkl operator
\be
{\cal D} = \partial_x + \nu x^{-1} (1-R) \lab{Dunkl} \ee
Consider e.g. the parametrization chosen in \re{L-1}, i.e.
\be
\lambda_n = (-1)^{n+1}(2n+\alpha+\beta+1)+\alpha+\beta+1, \quad \mu_n = (-1)^n (2n+\alpha)-\alpha \lab{lm1_-1} \ee 
Then 
\be
\nu = \alpha/2, \;\tilde \alpha=\alpha, \; \tilde \beta= \beta+2. \lab{til_-1} \ee
This means that 
\be
{\cal D} = \partial_x + \frac{\alpha} {2x} (1-R) \lab{Dunkl_1} \ee
and 
\be
\tilde \lambda_n = (-1)^{n+1}(2n+\alpha+\beta+3)+\alpha+\beta+3, \quad \tilde \mu_n = (-1)^n (2n+\alpha)-\alpha \lab{tlm_D} \ee
 The umbral classical property of the little -1 Jacobi polynomials with respect to the Dunkl operator   \re{Dunkl_1} was established in \cite{VZ_little}.

The more general problem is to find {\it all} admissible operators $\cal D$ and {\it all}  systems of umbral classical polynomials $P_n(x)$ satisfying property \re{tPD}. This problem is more complicated and remains open.

\section{Algebraic relations between the operators $L$ and $x$}
\setcounter{equation}{0}
The operator $L$ is low triangular and bidiagonal in the monomial basis as sen from \re{hyp-2-diag}. Introduce the operator $X$ which is multiplication by $x$.
Clearly
\be
X x^n = x^{n+1}, \lab{Xx} \ee 
i.e. the operator $X$ is the shift operator in the monomial basis.

Consider the operator 
\be
R_1 = X^2  L + L X^2 - \Omega XLX \lab{R_1} \ee
It is elementary to see that
\be
R_1 x^n =(\lambda_{n+2} + \lambda_n - \Omega \lambda_{n+1})x^{n+2} + (\mu_{n+2} + \mu_n - \Omega \mu_{n+1}) x^{n+1} \lab{R1xn} \ee
If $L$ is the admissible operator then its coefficients $\lambda_n, \mu_n$ satisfy relations \re{lin_lam}, \re{lin_nu}. In this case we can present \re{R1xn} in the algebraic form
\be
R_1 =-B_2 X^2 - C_2 X \lab{R1_XX} \ee
We thus have the 
\begin{pr}
The abstract "hypergeometric" operator $L$ is admissible if and only if there exists 3 parameters $\Omega, B_2,C_2$ such that the operator identity 
\be
 X^2  L + L X^2 - \Omega XLX +B_2 X^2 +C_2 X=0 \lab{crit_2} \ee
is valid on all polynomials.
\end{pr} 
Thus Proposition 1 and Proposition 2 give two equivalent necessary and sufficient conditions for the operator $L$ to be admissible.

On can introduce the "dual" operator
\be
R_2 = L^2 X + X L^2 - \Omega LXL \lab{R_2} \ee
It is directly verified that if the operator $L$ is admissible then relation
\be
R_2 = \alpha (LX + XL) +\beta X + \gamma L + \delta, \lab{R2_XL} \ee
holds, where $\alpha,\beta,\gamma, \delta$ are some constants.

Relations \re{R1_XX}- \re{R2_XL} can be considered as special case of generic Askey-Wilson algebra AW(3) introduced in \cite{ZheAW}, or equivalently, of the Askey-Wilson relations introduced by Terwilliger \cite{Ter}.

It is interesting to note that the monomial basis $x^n$ resembles the so-called split basis introduced in \cite{Ter}: the operator $L$ is low triangular while the operator $X$ is upper triangular in this basis. There is, however,  an essential difference with respect to \cite{Ter}. Indeed, in \cite{Ter} it was assumed that both operators $L$ and $X$ are finite-dimensional and diagonalizable. In our approach both operators $L,X$ are infinite-dimensional and moreover, the operator $X$ is non-diagonalizable.

\section{Conclusions}
\setcounter{equation}{0}
We have classified all admissible abstract "hypergeometric" operators $L$ having orthogonal polynomials $P_n(x)$ as eigensolutions. Apart from well known classical and q-classical polynomials there exists one more family of orthogonal polynomials (little -1 Jacobi polynomials) which are eigensolutions of the Dunkl type differential operator of the first order. In fact, the case corresponding to $\Omega  \ne \pm 2$ is generic (little q-Jacobi polynomials and their special and degenerate cases), while the cases corresponding to $\Omega = \pm 2$ can be obtained by corresponding limiting processes $q \to \pm 1$ (for details of the limit $q \to -1$ see, e.g. \cite{VZ_little}).

There are obvious generalizations of the above approach.

The first one is to consider the operators $L$ having more than 2 diagonals, say
\be
L x^n = \lambda_n x^n + \mu_n x^{n-1} + \nu_n x^{n-2} \lab{3-diag_L} \ee 
This operator preserves the space of polynomials and hence there exists a unique set of eigenpolynomials $P_n(x)$ satisfying the equation $L P_n(x) = \lambda_n P_n(x)$. It is natural to classify all the operators \re{3-diag_L} having orthogonal polynomials as eigensolutions. In contrast to "hypergeometric" operators with the property \re{hyp-2-diag} the expression for the polynomials $P_n(x)$ is not so simple. Nevertheless, we can consider the orthogonal polynomials $P_n(x)$ as having the "bispectrality" property. Recall that the orthogonal polynomials $P_n(x)$ are called the bispectral if they satisfy an additional eigenvalue problem \cite{Grun},\cite{SVZ}.
\be
L P_n(x) = \lambda_n P_n(x), \lab{bispec_L} \ee  
where the operator $L$ acts on the argument $x$ and can be either differential or difference operator of arbitrary order.

In \cite{Ismail} was suggested the bispectral problem \re{bispec_L} where 
\be
L= f(x) T + g(x) S + h(x) I. \lab{L_Ism} \ee
In \re{L_Ism} it is assumed that $S$ and $T$ are linear operators which map any polynomial of exact degree $n$ to polynomials of exact degrees $n-1$ and $n-2$ respectively, $I$ is the identical operator   and $f(x),g(x),h(x)$ are fixed polynomials. It can be showed \cite{Ismail} that all polynomials from the Askey scheme satisfy the bispectral property \re{bispec_L} with $L$ given by \re{L_Ism}. The concrete form of the operators $S,T$ depend on choice of the polynomials in the Askey scheme. For example, for the classical orthogonal polynomials $S=\partial_x, \: T=\partial_x^2$; for classical polynomials on the uniform grid $S =\Delta, \: T = \Delta \nabla$, where $\Delta f(x) = f(x+1)-f(x), \; \nabla f(x) = f(x)-f(x-1)$ etc. I.e. in all cases the operators $S,T$ belong to the class of either differential or difference operators.     

In our case the operator $L$  does not belong {\it a priori} to the class of differential or difference operators. Hence the bispectral problem \re{bispec_L} can be considered as a natural generalization of the classical bispectral problem.

Another generalization was mentioned in the previous section. Abstract "hypergeometric" polynomials satisfy "umbral" classical property   \re{tPD}. The inverse is not true: there are umbral classical polynomials beyond the scheme described in the present paper. It is expected that the family of umbral classical orthogonal polynomials is much wider than the family of classical or q-classical polynomials (including all polynomials from the Askey-Wilson scheme \cite{KLS}).

\bigskip\bigskip
{\Large\bf Acknowledgments}
\bigskip

\noindent The author is grateful to V.Genest, V. Spiridonov, S.Tsujimoto and L.Vinet for
discussion.

\newpage

\bb{99}


\bi{Al-Salam} W.A. Al-Salam, {\it Characterization theorems for orthogonal polynomials}, in: P. Nevai (Ed.), Orthogonal Polynomials: Theory and Practice, NATO ASI Series C: Mathematical and Physical Sciences, vol. {\bf 294}, Kluwer Academic Publishers, Dordrecht, pp. 1–-24.

\bi{AW} R.~Askey and J.~Wilson, {\it Some basic hypergeometric orthogonal
polynomials that generalize Jacobi polynomials}, Mem. Amer. Math. Soc. {\bf
54}, No. 319, (1985), 1-55.





\bi{Chi} T. Chihara, {\it An Introduction to Orthogonal
Polynomials}, Gordon and Breach, NY, 1978.

\bi{Duran} A.~J.~Dur\'an, {\it Orthogonal polynomials satisfying higher-order difference equations}, Constr.Approx.  {\bf 36} (2012), 459-–486.




\bi{Grun} F.Alberto Gr\"unbaum and L.Haine, {\it Bispectral Darboux transformations: an extension of the Krall polynomials}, IMRN, 1997, No. 8, 359--392.





\bi{Ismail} M.E.H.Ismail, {\it Classical and Quantum orthogonal polynomials in one variable}.
Encyclopedia of Mathematics and its Applications (No. 98), Cambridge, 2005.

\bi{Khol} A.N.Kholodov, {The umbral calculus and orthogonal polynomials}, Acta Appl.Mathem. {\bf 19} (1990), 1--54.

\bibitem{KLS} R. Koekoek, P.A. Lesky, and R.F. Swarttouw. {\it Hypergeometric orthogonal polynomials and their q-analogues}.
Springer, 1-st edition, 2010.





\bi{NSU} A.F. Nikiforov, S.K. Suslov, and V.B. Uvarov, {\em
Classical Orthogonal Polynomials of a Discrete Variable},
Springer, Berlin, 1991.




\bi{Roman} S.Roman, {\it The Theory of the Umbral Calculus. I }, J.Math.Anal.Appl. {\bf 87} (1982), 58--115.

\bi{SVZ} V.Spiridonov, L.Vinet and A.Zhedanov, {\it Bispectrality and Darboux transformation in the theory of orthogonal polynomials}, 111--122, in:  The Bispectral Problem, CRM Proceedings and Lecture Notes, v. 14, AMS, 1998.

\bi{SZ_root} V.Spiridonov, A.Zhedanov, {\it Zeros and orthogonality of the Askey-Wilson polynomials for $q$ a root of unity},  Duke Math. J. {\bf 89} (1997), 283--305.

\bi{Ter} P. Terwilliger, {\it Two linear transformations each tridiagonal with respect to an eigenbasis of the other}. Linear
algebra and its applications, {\bf 330} (2001), 149–-203.

\bi{VYZ} L.Vinet, O.Yermolayeva and A.Zhedanov, {\it A method to
study the Krall and q-Krall polynomials}, J.Comp.Appl.Math. {\bf
133} (2001) 647--½656.

\bi{VZ_Bochner} L.Vinet and A.Zhedanov, {\it Generalized Bochner
theorem: characterization of the Askey-½Wilson polynomials},
J.Comp.Appl.Math., {\bf 211} (2008) 45 -- 56.

\bi{VZ_little} L.Vinet and A.Zhedanov, {\it A 'missing' family of classical orthogonal polynomials}, J. Phys. A: Math. Theor. {\bf 44} (2011), 085201. arXiv:1011.1669


\bi{ZheAW}  A. Zhedanov, {\it Hidden symmetry of Askey-Wilson polynomials}. Theoretical and Mathematical Physics, {\bf 89} (1991) 1146–-1157.

\eb

\end{document}